\documentclass[11pt,fleqn,twoside]{article}
\usepackage{amsfonts,amssymb,latexsym,graphics}
\makeatletter
\newcommand{\prava}[1]{\small\it
\begin{flushleft}
Copyright \copyright \ 1999 by  #1
\end{flushleft}}

\newcommand{\name}[1]{\begin{flushleft}
                       \LARGE \bf #1
                       \end{flushleft}\vspace{-3mm}}

\newcommand{\Author}[1]{\begin{flushleft}
                       \it #1 \end{flushleft}}

\newcommand{\Adress}[1]{\begin{flushleft}
                       \it #1 \end{flushleft}}

\newcommand{\Date}[1]{\begin{flushleft}
                      \small  \it #1 \end{flushleft}}

\newcommand{\ehkol}{Author \ name}
\newcommand{\ohkol}{Article \ name}
\renewcommand{\@evenhead}{
\hspace*{-3pt}\raisebox{-15pt}[\headheight][0pt]{\vbox{\hbox to \textwidth 
{\thepage \hfil \ehkol}\vskip4pt \hrule}}}
\renewcommand{\@oddhead}{
\hspace*{-3pt}\raisebox{-15pt}[\headheight][0pt]{\vbox{\hbox to \textwidth 
{\ohkol \hfil \thepage}\vskip4pt\hrule}}}
\renewcommand{\@evenfoot}{}
\renewcommand{\@oddfoot}{}

\newcommand{\be}{\begin{equation}}
\newcommand{\ee}{\end{equation}}
\newcommand{\ba}{\hspace*{-5pt}\begin{array}}
\newcommand{\ea}{\end{array}}

\newcommand{\ds}{\displaystyle}
\makeatother

\begin{document}
\setcounter{page}{35}
\thispagestyle{empty}

\renewcommand{\ehkol}{V.A. Trotsenko}
\renewcommand{\ohkol}{Variational Methods for Solving Nonlinear Boundary 
Problems of Statics}

\begin{flushleft}
\footnotesize \sf
Journal of Nonlinear Mathematical Physics \qquad 1999, V.6, N~1,
\pageref{trotsenko-fp}--\pageref{trotsenko-lp}. 
\hfill {\sc Article}
\end{flushleft}

\vspace{-5mm}

\renewcommand{\footnoterule}{}
{\renewcommand{\thefootnote}{} \footnote{\prava{V.A. Trotsenko}}}

\name{Variational Methods for Solving \\
Nonlinear Boundary Problems of\\
Statics of Hyper-Elastic Membranes}\label{trotsenko-fp}

\Author{V.A.~TROTSENKO}

\Adress{Institute of Mathematics of the National Academy of Sciences of
Ukraine, \\
3 Tereshchenkivska Str., 252601 Kyiv-4, Ukraine}

\Date{Received June 10, 1998; Accepted July 13, 1998}

\begin{abstract}
\noindent
A number of important results of studying large deformations of
hyper-elastic shells are obtained using discrete
methods of mathematical physics \cite{trotsenko:1:}--\cite{trotsenko:6:}.
In the present paper, using  the variational method for
solving nonlinear boundary
problems of statics of hyper-elastic membranes  under the regular
hydrostatic load, we investigate peculiarities of deformation of a circular
membrane whose mechanical characteristics are described by the
Bidermann-type elastic potential.
We develop an algorithm for solving a singular perturbation of
nonlinear problem for the case of membrane loaded by heavy liquid.
This algorithm enables us to  obtain approximate solutions both in the
presence of boundary layer and without it.
The class of  admissible functions, on which the variational method is
realized, is chosen  with account of the  structure of
formal asymptotic expansion of solutions of the corresponding
linearized equations that have singularities
in a small parameter at higher derivatives and in the
independent variable. We give examples of calculations
that  illustrate possibilities of the method suggested for solving
the problem  under consideration.
\end{abstract}

\setcounter{equation}{0}
\renewcommand{\theequation}{\arabic{section}.\arabic{equation}}

\section{Def\/inition of the strained deformed state of a circular
membrane under the regular hydrostatic load}

Let us consider a circular membrane of the radius $R_0$ of
incompressible, isotropic and hyperelastic material having the
small constant
width $h_0$. Let a hydrostatic load $Q$ be applied to the
membrane with an inf\/lexibly f\/ixed contour. To describe
the geometry
of the deformed membrane, let us introduce the cylindrical
coordinate system $Oz\eta r$, the axis $Oz$ of which coincides
with the  symmetry axis of the membrane.
It follows from the axial symmetry of
the problem that principal deformation directions at any point
will coincide with meridians, parallels and normals to the
deformed surface and all the parameters of the strained deformed
state will be functions of  initial distances $s$ of
points of the membrane from its symmetry axis only. Denote the
principal degrees of lengthenings in these directions by
$\lambda_1$, $\lambda_2$ and $\lambda_3$, respectively. Then,
according to the nonlinear theory of elastic membrane \cite{trotsenko:1:},
interior stresses $T_1$ and $T_2$ of the deformed shell in
the direction of its meridian and parallel can be determined
from the relations
\begin{equation} \label{trotsenko:1.1}
\ba{l}
\displaystyle
T_i=2h_0\lambda_3\left(\lambda_i^2-\lambda_3^2\right)
\left(\frac{\partial W}{\partial I_1}+\lambda_{3-i}^2
\frac{\partial W}{\partial I_2}\right)
\qquad (i=1,2),
\vspace{3mm} \\
\displaystyle
\lambda_1=\sqrt{\left(\frac{dz}{ds}\right)^2+\left(\frac{dr}{ds}\right)^2},
\qquad
\lambda_2=\frac{r(s)}{s},
\qquad
\lambda_3=\frac{1}{\lambda_1\lambda_2}=\frac{h(s)}{h_0},
\vspace{3mm} \\
\displaystyle
I_1=\lambda_1^2+\lambda_2^2+\lambda_3^2,
\qquad
I_2=\lambda_1^{-2}+\lambda_2^{-2}+\lambda_3^{-2},
\quad
W=W(I_1,I_2).
\ea
\end{equation}
Here $z(s)$ and $r(s)$ are axial and radial components of the
radius of the generatix vector of the deformed shell, $W$ is
a deformation energy function for the material of the
shell, $I_1$ and $I_2$ are deformation invariants, and $h(s)$ is the
width of the membrane in the deformed state.

Let us choose the deformation energy function $W(I_1,I_2)$ in
the form of its four-term approximation proposed by
V.L.~Biderman \cite{trotsenko:7:}, i.e., 
\begin{equation} \label{trotsenko:1.2}
W=C_1(I_1-3)+C_2(I_2-3)+C_3(I_1-3)^2+C_4(I_1-3)^3,
\end{equation}
where $C_i$ $(i=1,\ldots,4)$ are physical constants
that are  determined experimentally.

By using function (\ref{trotsenko:1.2}) one can suf\/f\/iciently precisely
describe the law of deforming for a number of types of elastic
materials up to their break.

If the origin of the coordinate system coincides with a pole of
the deformed membrane, then the hydrostatic pressure acting on
the shell, is of the form
\begin{equation} \label{trotsenko:1.3}
Q=C^*-Dz,
\qquad
C^*=P^*-P^0,
\qquad
D=\rho g.
\end{equation}
Here $\rho$ is the density of liquid, $g$ is the acceleration of
gravity, $P^*$ is a constant component of pressure in liquid,
and
$P^0$ is the pressure of gas onto the shell.

Let us introduce the following dimensionless parameters:
\[
\ba{l}
\{s^*,r^*,z^*\}=\{s,r,z\}/R_0,
\qquad
T_i^*=T_i/(2C_1h_0),
\\[2mm]
W^*=W/C_1,
\qquad
Q^*=QR_0/(2C_1h_0).
\ea
\]

In what follows, we shall use these parameters
omitting asterisks for the sake of simplicity.

Movements corresponding to the static equilibrium state are
singled out from all admissible movements by the fact that
 stationary values of the functional
\cite{trotsenko:8:}
\begin{equation} \label{trotsenko:1.4}
I=\int\limits_0^1\left[W(I_1,I_2)s-Qr^2\frac{dz}{ds}\right]ds
\end{equation}
correspond to them.

Having placed the origin of the coordinate system $Oz\eta r$ at the
center of the undeformed membrane, we present the
hydrostatic load as
\[
Q=C-Dz,
\]
where dif\/ference between the  constants $C$ and
$C^*$ in formula (\ref{trotsenko:1.3}) is due to  the shift of the
origin of the coordinate system. In this case, we should seek a
stationary value of functional (\ref{trotsenko:1.4}) within the set of
functions satisfying the boundary conditions
\begin{equation} \label{trotsenko:1.5}
\frac{dz}{ds}\biggr|_{s=0}=0,
\qquad
z(1)=0, \qquad r(0)=0, \qquad r(1)=1.
\end{equation}

Utilizing usual tools of the calculus of variations, one can
show that the Euler equations for functional (\ref{trotsenko:1.4}) are:
\begin{equation} \label{trotsenko:1.6}
\frac{d}{ds}(rT_1)=T_2\frac{dr}{ds},
\qquad
k_1T_1+k_2T_2=Q.
\end{equation}
Here $k_1$ and $k_2$ are principle curvatures of the deformed
median surface of the membrane calculated  by the formulae
\[
k_1=\left(\frac{d^2r}{ds^2}\frac{dz}{ds}-
\frac{dr}{ds}\frac{d^2z}{ds^2}\right)\lambda_1^{-3},
\qquad
k_2=-(r\lambda_1)^{-1}\frac{dz}{ds}.
\]

Following Rietz to
f\/ind an extremum of functional (\ref{trotsenko:1.4})
we present functions $z(s)$ and $r(s)$ in the form of the
following expansions:
\begin{equation} \label{trotsenko:1.7}
z(s)=\sum\limits_{k=1}^mx_ku_k(s),
\qquad
r(s)=s+\sum\limits_{k=1}^mx_{k+m}v_k(s),
\end{equation}
where in view of (\ref{trotsenko:1.5}), the  coordinate
functions $u_k(s)$ and $v_k(s)$  have to obey the conditions
\begin{equation} \label{trotsenko:1.8}
u_k^\prime(0)=u_k(1)=v_k(0)=v_k(1)=0.
\end{equation}
Then constants $x_k$ $(k=1,2,\ldots,2m)$ forming $2m$-component
vector $\vec x$ will be determined from  the
nonlinear algebraic system
\begin{equation} \label{trotsenko:1.9}
\vec g(\vec x)=0.
\end{equation}
In this case, components of the vector-function $\vec g$ are of
the form
\begin{equation} \label{trotsenko:1.10}
\ba{l}
\displaystyle
g_i=\int\limits_0^1\left[U(\lambda_1,\lambda_2)\frac{dz}{ds}
\frac{du_i}{ds}-Q\lambda_2\frac{dr}{ds}u_i\right]s\hspace{1.5pt}ds,
\vspace{3mm} \\
\displaystyle
g_{i+m}=\int\limits_0^1\left[U(\lambda_1,\lambda_2)\frac{dr}{ds}
\frac{dv_k}{ds}+\left(U(\lambda_2,\lambda_1)
\frac{\lambda_2}{s}+Q\lambda_2\frac{dz}{ds}\right)v_i\right]
s\hspace{1.5pt}ds.
\ea
\end{equation}
Here
\[
U(\lambda_1,\lambda_2)=\left(1-\lambda_1^{-4}\lambda_2^{-2}\right)
\left(\frac{\partial W}{\partial I_1}+\lambda_2^2
\frac{\partial W}{\partial I_2}\right).
\]

We  solve algebraic system (\ref{trotsenko:1.9}) by using the Newton
method according to which correction of an approximate solution
at each iteration step is performed by the scheme
\begin{equation} \label{trotsenko:1.11}
\vec x^{(k+1)}=\vec x^{(k)}-H^{-1}\left(\vec x^{(k)}\right)
g\left(\vec x^{(k)}\right)
\quad
(k=1,2,\ldots), 
\end{equation}
where $H(\vec x)$ is the symmetric Jacobi matrix of functions
$g_1$, $g_2,\ldots,g_{2m}$ of variables $x_1$,
$x_2,\ldots,x_{2m}$, whose elements  are calculated as follows
\[
h_{ij}=\int\limits_0^1\left\{\left[
\frac{\partial U(\lambda_1,\lambda_2)}{\lambda_1\partial\lambda_1} 
\left(\frac{dz}{ds}\right)^2+U(\lambda_1,\lambda_2)\right]
\frac{du_i}{ds}\frac{du_j}{ds}+D\lambda_2\frac{dr}{ds}u_iu_j\right\}
s\hspace{1.5pt}ds,
\]
\[
h_{i+m,j+m}=\int\limits_0^1\left\{\left[
\frac{\partial U(\lambda_1,\lambda_2)}{\lambda_1\partial\lambda_1} 
\left(\frac{dz}{ds}\right)^2+U(\lambda_1,\lambda_2)\right]
\frac{dv_i}{ds}\frac{dv_j}{ds}  \right.
\]
\[
\qquad
\left. +\frac 1s\frac{dr}{ds}
\frac{\partial U(\lambda_1,\lambda_2)}{\partial\lambda_2} 
\frac{d}{ds}(v_iv_j)+
\left[\frac{\partial U(\lambda_2,\lambda_1)}
{\lambda_2\partial\lambda_2}+ U(\lambda_2,\lambda_1)+
Qs\frac{dz}{ds}\right]\frac{v_iv_j}{s^2}\right\}s\hspace{1.5pt}ds,
\] \[
\qquad
(j=\overline{1,m}; \ i\geq j),
\] \[
h_{i,m+j}=\int\limits_0^1\left\{
\frac{\partial U(\lambda_1,\lambda_2)}{\lambda_1\partial\lambda_1} 
\frac{dz}{ds}\frac{dr}{ds}\frac{du_i}{ds}\frac{dv_j}{ds}\right.
\] \[
\qquad
\left. +\left[
\frac{\partial U(\lambda_1,\lambda_2)}{s\partial\lambda_2}
\frac{dz}{ds}+Q\lambda_2\right]
\frac{du_i}{ds}v_j-D\lambda_2\frac{dz}{ds}u_iv_j\right\}
s\hspace{1.5pt}ds
\qquad
(i,j=1,\ldots,m).
\]

So starting  from the potential energy stationarity
principle, we reduce the problem of determining f\/inite deformations of the
circular membrane  to calculating integrals
forming algebraic system and elements of the Jacobi matrix
with subsequent solving problems of linear algebra (\ref{trotsenko:1.11}) at each
step of the Newton iteration process. The algorithm may be
ef\/fective if one chooses properly  systems of coordinate functions
$\{u_k(s)\}$ and $\{v_k(s)\}$ that would allow one to get
solutions with high precision, provided the dimension of the
nonlinear algebraic system is  not large.

To this end, let us expand  solutions to be found in a Taylor
series in a neighbourhood of the point $s=0$,  using
sequential dif\/ferentiation of equilibrium equations and
geometric relationships of the shell in order to get coef\/f\/icients of
the expansion. With regard to symmetry conditions at the pole of
the deformed shell
\[
\lambda_1=\lambda_2=\lambda,
\qquad
T_1=T_2=T,
\qquad
k_1=k_2=k
\qquad
(s=0),
\]
we can show that solutions boundered under $s\to 0$ have the
following structure:
\begin{equation} \label{trotsenko:1.12}
\begin{array}{l}
z(s)=a_1+a_2s^2+a_3s^4+\cdots,
\vspace{2mm} \\
r(s)=b_1s+b_2s^3+b_3s^5+\cdots.
\end{array}
\end{equation}

Taking into account expansions (\ref{trotsenko:1.12}), sequences of
coordinate functions $\{u_k(s)\}$ and $\{v_k(s)\}$ satisfying
boundary conditions (\ref{trotsenko:1.8})  take the form
\begin{equation} \label{trotsenko:1.13}
u_k(s)=(s^2-1)s^{2k-2},
\qquad
v_k(s)=(s^2-1)s^{2k-1}
\qquad
(k=1,\ldots,m).
\end{equation}

Let us present some results of calculations according to the
algorithm proposed. To illustrate convergence of the Rietz
method,  we write out values of functions~$z(s)$ and
$r(s)$ and their f\/irst two derivatives, depending on the number $m$
of approximations in expansions (\ref{trotsenko:1.7}) at the point $s=0.2$
for $C=1.7$, $D=0$ in Table~1. Ratios of elastic constants in Biderman's
potential are chosen as follows:
\[
\Gamma_1=0.02, \qquad \Gamma_2=-0.015, \qquad \Gamma_3=0.00025.
\]
In the last column of the table we give the relative error of
the obtained approximation to the solution of
the second equilibrium equation
\begin{equation} \label{trotsenko:1.14}
\Delta=|k_1T_1+k_2T_2-Q|/|C|.
\end{equation}

\begin{center}
{\bf Table 1}\\[3mm]
\begin{tabular}{cccccccc}
\hline
\multicolumn{1}{c|}{$m$\rule[-7pt]{0pt}{20pt}} & \multicolumn{1}{c|}{$z$}
& \multicolumn{1}{c|}{$r$} & \multicolumn{1}{c|}{$-z'$} &
\multicolumn{1}{c|}{$r'$} & \multicolumn{1}{c|}{$-z''$} &
\multicolumn{1}{c|}{$-r''$} & \multicolumn{1}{c}{$\Delta$} \\
\hline
\rule{0pt}{13pt}%
1 & 0.7016 & 0.2865 & 0.2923 & 1.3966 & 1.4617 & 0.5408 &
$2\cdot 10^{-1}$ \\
2 & 0.7824 & 0.3026 & 0.3975 & 1.4644 & 1.9279 & 0.7252 &
$6\cdot 10^{-2}$ \\
3 & 0.7913 & 0.3064 & 0.4255 & 1.4753 & 2.0248 & 0.8345 &
$1\cdot 10^{-2}$ \\
4 & 0.7926 & 0.3068 & 0.4350 & 1.4759 & 2.0421 & 0.8560 &
$1\cdot 10^{-3}$ \\
5 & 0.7926 & 0.3069 & 0.4362 & 1.4758 & 2.0415 & 0.8591 &
$6\cdot 10^{-4}$ \\
6 & 0.7926 & 0.3069 & 0.4362 & 1.4757 & 2.0406 & 0.8596 &
$2\cdot 10^{-5}$ \\
\hline
\end{tabular}
\end{center}


Solutions possess the similar convergence at other points of the
interval of integration for these equations. Lengthening of the
generatrix of the deformed shell af\/fects weakly on the degree of
convergence of the variational method.  For a chosen system of coordinate
functions, solutions have the uniform convergence as well as their
f\/irst two derivatives for suf\/f\/iciently wide range of parameters.

\begin{figure}[t]
\centerline{\scalebox{1}{\includegraphics{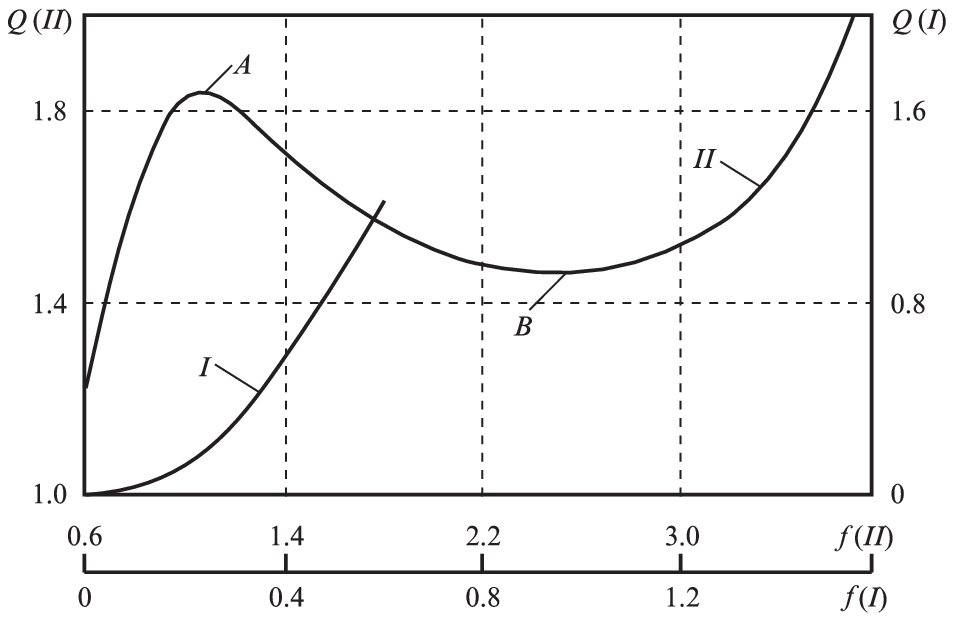}}}
\begin{center}
Figure 1
\end{center}

\vspace{-7mm}

\end{figure}

Graphic dependence of the sag of the shell at its pole $f=z(0)$,
upon the applied uniform pressure $(D=0)$ is shown in Figure\,1. The
characteristic indicates  the possibility of existence of several equilibrium
forms for the same value of the load. It is clearly seen that the
diagram "load-sag" for the material under consideration has two
extreme points $A$ and $B$. If  the
parameter of the load increases monotonely and reaches the
value of the upper critical load,
then there will be an abrupt passage of the system into a new
equilibrium state, after which sags of the shell will again
increase smoothly. The monotone decrease of the load $C$
(unloading stage) leads to the backward distortion of the
shell to the initial equilibrium form at lower critical load.
Equilibrium states described by a falling part of the curve from
the point $A$ to the point $B$ are unstable and they are not
realized at all in the way of loading under consideration.  

In conclusion, we consider  to the main principles of
choosing an initial approximation of the Newton iteration
process that allows us to get solutions of nonlinear algebraic
system (\ref{trotsenko:1.9}). If solutions of these equations are far from
singular points $A$ and~$B$, then one can construct an initial
approximation, taking advantage of rapid convergence of
solutions. As we can see from Table 1, the f\/irst Rietz
approximation is a good solution to the initial problem. In this
case, we have to solve a system of two nonlinear algebraic
equations. The initial approximation is determined from
physically noncontradictory values of $x_1$ and~$x_2$ (to the
positive load there corresponds the positive sag of the
shell). Then we use the obtained solution to construct $P$-th
approximation, setting the f\/irst and $(P+1)$-th components of
the vector $\vec x$ equal to obtained values of $x_1$ and $x_2$,
and putting other components equal to zero.
The initial approximation constructed in
this way improved by employing the Newton method
with 4--5 iterations. In the case, when it is necessary to obtain
a series of equilibrium states, when choosing initial values of
the vector $\vec x$ one can use dif\/ferent extrapolational
formulae to prolong solutions in the parameter~$C$. However,
when implementing the procedure of prolonging a solution with
approximation to critical points $A$ and $B$, the convergence of
the Newton process sharply comes down, and then the divergence
appears clearly. To go through these equilibrium
states of the shell one should change prolongation
parameters. The main point of this approach is that starting
from some moment, the parameter of the load is supposed to be
unknown, and as a given parameter we choose the value of the sag
of the shell at the pole. Artif\/icial introduction of a new
parameter provides that the uniqueness condition for a solution
in a neighbourhood of a f\/ixed parameter is fulf\/illed, and,
therefore, Jacobian is nondegenerate, which is a necessary
condition for  the Newton method to be applicable.

The algorithm proposed is also ef\/fective for other rotation
shells in the nondeformed state, the generatrix of which
intersects the symmetry axis at the right angle and principle
curvatures at its pole are equal.

\setcounter{equation}{0}

\section{Adaptive Rietz method under conditions of singular
perturbance of the initial problem}

Let us consider the case of loading a membrane by a hydrostatic
load for $C^*\ll 1$ and $D\gg 1$ (heavy liquid). For such
parameters of the load, stresses in the median surface of the
deformed shell will be small, i.e., $\{T_1,T_2\}\ll 1$, and
the shell will take a nearly plane form with the exception of
a narrow domain near the f\/ixed contour of the membrane. The
increase of the constant component of the load $C^*$ leads to the
increase of stresses in the membrane, that begin to play the
prevalent role in forming its deformed surface. This is
accompanied by increasing the width of a domain with rapidly
changing functions of movements and stresses, and turning a plane
part of the surface of the shell into a curvilinear surface.

After introducing a small parameter $\mu=1/D$, equilibrium
equations (\ref{trotsenko:1.6}) in the coordinate system associated with the
center of the nondeformed membrane can be rewritten as follows:

\begin{equation} \label{trotsenko:2.1}
\begin{array}{l}
L_1(s,r,z',r',z'',r'')=0,
\vspace{2mm} \\
\mu\left[L_2(s,r,z',r',z'',r'')-C^*\right]=l-z,
\end{array}
\end{equation}
where 
\[
L_1=\frac{dT_1}{ds}+\frac 1r\frac{dr}{ds}(T_1-T_2),
\qquad
L_2=k_1T_1+k_2T_2,
\qquad
l=z(0).
\]

In the case under consideration, the presence of a small
parameter at higher derivatives transfers the initial regular
problem into a class of nonlinear singular perturbed problems.
In this case, solutions will  have domains both
with fast and slow variability. As a result, when
approximating solutions by means of polynomial basis, the
convergence of successive approximations in the Rietz method
sharply comes down, and the increase of the number of coordinate
functions in expansions (\ref{trotsenko:1.7}) leads to a loss of
calculation stability before the required precision is
achieved. The situation is also aggravate by that the
singularity degree of the problem depends on a solution itself
that is a priori unknown. Thus, there arise the principal problem  of
extending a class of admissible functions by functions that allow
getting
within the framework of a unique algorithm  approximate solutions
with the equally high precision both in the presence of a
narrow area of great changes in solutions and without it.

The problem indicated plays a key role in realizing the method of
variations when solving singularly perturbed boundary problems.
It is one of suf\/fuciently dif\/f\/icult problems of mathematical
physics and theory of function approximation. Apparently for
this reason, variantional methods (unlike f\/inite dif\/ference
methods \cite{trotsenko:9:}) are not developed properly for
solving even
linear boundary problems with a small parameter  the higher
derivative. 

One of the most important merits of the method of variations
is that when constructing approximate solutions,
one can use  an information about  problem, that can be
obtained by means of preliminary analysis of required solutions.
On the other hand, methods of a small parameter  may be
also used to  determine a structure
of a solution and characterize its degeneration with the
parameter at the higher derivative tending to zero. With this
objective  in view, let us carry out an asymptotic analysis of the
initial equations.

Searching for a direct expansion of solutions to system (\ref{trotsenko:2.1})
in the form
\begin{equation} \label{trotsenko:2.2}
z(s)=\sum\limits_{k=0}^\infty\mu^kz^{(k)}(s),
\qquad
r(s)=\sum\limits_{k=0}^\infty\mu^kr^{(k)}(s),
\end{equation}
let us expand operators $L_i$ $(i=1,2)$ in a Taylor series in the
parameter $\mu$
\begin{equation} \label{trotsenko:2.3}
L_i=(L_i)_{\mu=0}+\frac{\mu}{1!}\left(\frac{dL_i}{d\mu}\right)_{\mu=0}+
\cdots+\frac{\mu^{(n)}}{n!}\left(\frac{d^nL_i}{d\mu^n}\right)_{\mu=0}+
\cdots.
\end{equation}
Having substituted (\ref{trotsenko:2.2}) and (\ref{trotsenko:2.3})
into (\ref{trotsenko:2.1})
and having equated coef\/f\/icients of expansions of the same degrees of
$\mu$, we come to
\begin{equation} \label{trotsenko:2.4}
\ba{l}
\displaystyle
(L_1)_{\mu=0}=0,
\qquad z^{(0)}(s)=l,
\vspace{3mm} \\
\displaystyle
\left(\frac{dL_1}{d\mu}\right)_{\mu=0}=0,
\qquad
(L_2)_{\mu=0}-C^*=-z^{(1)}(s),
\vspace{3mm} \\
\displaystyle
\left(\frac{d^nL_1}{d\mu^n}\right)_{\mu=0}=0,
\qquad
\frac{1}{(n-1)!}\left(\frac{d^{(n-1)}L_2}{d\mu^{(n-1)}}\right)_{\mu=0}=
-z^{(n)}(s) 
\qquad (n=2,3,\ldots).
\ea
\end{equation}

We can obtain a bounded solution of the nonlinear equation
$(L_1)_{\mu=0}=0$ with respect to the function $r^{(0)}(s)$ by
expanding it in a Taylor series in a neighbourhood of the point~$s=0$
with regard to the symmetry at a pole of the deformed
shell. In this case, we can show that
\begin{equation} \label{trotsenko:2.5}
r^{(0)}(s)=a_1s,
\qquad
z^{(0)}(s)=l,
\end{equation}
where $a_1$ is an arbitrary constant.

Furthermore, with regard for expressions (\ref{trotsenko:2.5}), linear
equations of the f\/irst approximation allow us to get expressions
for functions  $r^{(1)}(s)$ and $z^{(1)}(s)$:
\[
r^{(1)}(s)=a_2s,
\qquad
z^{(1)}(s)=C^*.
\]

Having performed several transformations
that are simple but suf\/f\/iciently awkward, we can show that
\[
r^{(n)}(s)=a_ns,
\qquad
z^{(n)}(s)=0 
\qquad
(n=2,3,\ldots),
\]
where $a_n$ are arbitrary constants.

Collecting the results obtained, we come to the general form of a
solution bounded for $s=0$ as a result of direct expansion
(\ref{trotsenko:2.2}):
\begin{equation} \label{trotsenko:2.6}
r(s)=A(\mu)s,
\qquad
z(s)=C/D,
\qquad
C=C^*+lD,
\end{equation}
where $A(\mu)$ is a function of $\mu$; $C$ is the constant
component of the hydrostatic pressure, that dif\/fers from $C^*$
at the expense of that the origin of the coordinate system
passes into the center of the nondeformed membrane.

Proceeding from equations (\ref{trotsenko:2.1}), it is not possible to
establish the structure of solutions in the boundary area, since
equations of the f\/irst approximation are nonlinear and they are
not integrated in the explicit form. For this reason, in what
follows, using the results  of the asymptotic analysis of the linearized
equilibrium equations with respect to some equilibrium state
that is described by an analytic function of the form (\ref{trotsenko:1.12}),
we propose to establish the structure of their solutions and
after that to endow solutions of the nonlinear problem with
these properties. 

Suppose that the membrane is in the deformed state as a result
of applying the hydrostatic load $Q=C^*-Dz$. Let the shell
have passed into a new state at the expense of changing the
constant pressure component by a small value $\delta C^*$.
Denote by $u$ and $w$ projections of a movement of the shell
from this load onto directions of its generatrix and outer
normal, respectively. We shall characterize the main strained
deformed state by stresses $T_1$ and $T_2$, and principal
degrees of lengthenings $\lambda_1$, $\lambda_2$ and
$\lambda_3$. Having performed linearization of elasticity
relations (\ref{trotsenko:1.1}), we can establish the following
relation between additional stresses
$\delta T_1$ and $\delta T_2$ and deformations of the shell:
\begin{equation} \label{trotsenko:2.7}
\delta T_1=c_{11}\varepsilon_1+c_{12}\varepsilon_2,
\qquad
\delta T_2=c_{21}\varepsilon_1+c_{22}\varepsilon_2,
\end{equation}
where
\[
\varepsilon_1=\frac{1}{\lambda_1}\frac{du}{ds}+k_1w,
\qquad
\varepsilon_2=\frac{1}{r\lambda_1}\frac{dr}{ds}u+k_2w;
\]
\[
c_{11}=f_1(\lambda_1,\lambda_2),
\qquad
c_{12}=f_2(\lambda_1,\lambda_2),
\qquad
c_{21}=f_2(\lambda_2,\lambda_1),
\qquad
c_{22}=f_1(\lambda_2,\lambda_1);
\]
\[
f_1(\lambda_1,\lambda_2)=\left(\lambda_1^2\lambda_3+3\lambda_3^3\right)
\frac{\partial W}{\partial I_1}+
\left(\lambda_1\lambda_2+3\lambda_2^2\lambda_3^3\right)
\frac{\partial W}{\partial I_2}
\]
\[
\qquad
+2\lambda_3\left(\lambda_1^2-\lambda_3^2\right)\left(A_{11}+
2A_{12}\lambda_2^2+A_{22}\lambda_2^4\right),
\]
\[
f_2(\lambda_1,\lambda_2)=\left(3\lambda_3^3-\lambda_1^2\lambda_3\right)
\frac{\partial W}{\partial I_1}+
\left(\lambda_1\lambda_2+\lambda_2^2\lambda_3^3\right)
\frac{\partial W}{\partial I_2}
\]
\[
\qquad
+2\lambda_3\left(\lambda_1^2-\lambda_3^2\right)
\left(\lambda_2^2-\lambda_3^2\right)\left(A_{11}+
A_{12}\left(\lambda_2^2+\lambda_2^2\right)+
A_{22}\lambda_1^2\lambda_2^2\right),
\]
\[
A_{ik}=\frac{\partial^2W}{\partial I_i\partial I_k} 
\qquad
(i,k=1,2).
\]

Linearized equilibrium equations in movements, describing the
perturbed state of the shell are of the form
\begin{equation} \label{trotsenko:2.8}
\ba{l}
\displaystyle
L_{11}(u)+L_{12}(w)=0,
\vspace{3mm} \\
\displaystyle
\varepsilon^2\left[L_{21}(u)+L_{22}(w)\right]+
r\left(\frac{dr}{ds}w+\frac{dz}{ds}u\right)=
\left(\varepsilon^2\delta C^*+\Delta l\right)r\lambda_1,
\ea
\end{equation}
where
\[
\ba{l}
\ds L_{11}(u)=-\frac{d}{ds}\left(\alpha_1\frac{du}{ds}\right)+\alpha_2u,
\qquad
L_{12}(w)=\alpha_3\frac{dw}{ds}+\alpha_4w,
\vspace{3mm} \\
\ds L_{22}(w)=-\frac{d}{ds}\left(\beta_1\frac{dw}{ds}\right)+\beta_2w,
\qquad
L_{21}(u)=\beta_3\frac{du}{ds}+\beta_4u,
\vspace{3mm} \\
\ds \alpha_1=\frac{c_{11}r}{\lambda_1},
\qquad
\alpha_2=\alpha_2^{(1)}+\frac{d\alpha_2^{(2)}}{ds},
\qquad
\alpha_2^{(1)}=\frac{c_{22}}{r\lambda_1}\left(\frac{dr}{ds}\right)^2-
r\lambda_1k_1k_2T_2,
\vspace{3mm} \\
\ds \alpha_2^{(2)}=-\frac{c_{12}+T_1}{\lambda_1}\frac{dr}{ds},
\qquad
\alpha_3=-r(c_{11}k_1+c_{21}k_2),
\qquad
\alpha_4=\alpha_4^{(1)}+\frac{d\alpha_4^{(2)}}{ds},
\vspace{3mm} \\
\ds \alpha_4^{(1)}=\bigl[(c_{21}+T_2)k_1+c_{22}k_2\bigr]\frac{dr}{ds},
\qquad
\alpha_4^{(2)}=-r\bigl[c_{11}k_1+(c_{12}+T_1)k_2\bigr],
\vspace{3mm} \\
\ds \beta_1=\frac{rT_1}{\lambda_1},
\qquad
\beta_2=r\lambda_1\left[(c_{11}-T_1)k_1^2+(c_{12}+c_{21})k_1k_2+
(c_{22}-T_2)k_2^2\right],
\vspace{3mm} \\
\ds \beta_3=-\alpha_3,
\qquad
\beta_4=\beta_4^{(1)}+\frac{d\beta_4^{(2)}}{ds},
\qquad
\beta_4^{(1)}=\bigl[c_{12}k_1+(c_{22}-T_2)k_2\bigr]\frac{dr}{ds},
\vspace{3mm} \\
\ds \beta_4^{(2)}=rT_1k_1,
\qquad
\varepsilon^2=\mu, 
\qquad
\Delta l=w(0).
\ea
\]

One can show that the structure of a direct expansion of
solutions to equations (\ref{trotsenko:2.8}) in the parameter
$\varepsilon$ corresponds to  structure (\ref{trotsenko:2.6})
established for the system of initial nonlinear equations. The
subject of further research is to construct integrals with high
variability of a system of linear equations with variable
coef\/f\/icients (\ref{trotsenko:2.8}). Since a direct expansion is an
approximation to some solution of initial equations (\ref{trotsenko:2.8}),
we shall seek boundary layer part of an asymptotics as an
approximation to a solution of the relative homogeneous system
of equations. In so doing, we shall take into account that
functions~$r(s)$ and $z(s)$ characterizing the geometry of the
shell in the main deformed state are analytic functions in a
certain domain and according to (\ref{trotsenko:1.12}) they are expanded
in power series in even and odd powers of $s$, respectively. On
the basis of this fact, we can rewrite the initial homogeneous
system of equations as follows:
\begin{equation} \label{trotsenko:2.9}
\ba{l}
\displaystyle
s^2\frac{d^2u}{ds^2}+sp_{11}\frac{du}{ds}+
sp_{12}\frac{dw}{ds}+q_{11}u+q_{12}w=0,
\vspace{3mm} \\
\displaystyle
\varepsilon^2
s^2\frac{d^2w}{ds^2}+\varepsilon^2\left(sp_{21}\frac{du}{ds}+
sp_{22}\frac{dw}{ds}\right)+q_{21}u+q_{22}w=0,
\ea
\end{equation}
where $p_{ij}(s)$ and $q_{ij}(s)$ are analytic functions for 
$s\in [0,1]$.

Proceeding from the form of system of equations (\ref{trotsenko:2.9}), we
 conclude that this system has singularities both in the
parameter and unknown variable. The presence of a small
parameter at the higher derivative results in existence of a
narrow area in a neighbourhood of the point $s=1$, in which
solutions of the system have large gradients. The point $s=0$, in
its turn, is a regular point for equations (\ref{trotsenko:2.9}), that
assumes the relative asymptotics of solutions when $s\to 0$.
Since our f\/inal goal is to construct an algorithm for solving a
problem with the uniform convergence in the parameter
$\varepsilon$, in what follows it is necessary to take into
account both of these facts.

We have  not succeededin constructing solutions  a singular
perturbed equation, with a regular singular point, that belongs
to the class of equations under
consideration containing an
exponent factor \cite{trotsenko:10:}. Simplif\/ied equations obtained from
the initial ones while only the f\/irst terms of expansion of
their coef\/f\/icients remain unchanged, may serve as a
guiding line for an asymptotic representation of solutions to
system (\ref{trotsenko:2.9}). In this case, we can establish that a
regular solution of these equations is expressed by Bessel
functions of order zero and their f\/irst derivatives. According
to what has been said, we  seek solutions of homogeneous
system of equations (\ref{trotsenko:2.8}) with a large variability index
in the following form:
\begin{equation} \label{trotsenko:2.10}
\ba{l}
u=u_1(s,\varepsilon)J_0(x)+u_2(s,\varepsilon)J_0^\prime(x),
\qquad
w=w_1(s,\varepsilon)J_0(x)+w_2(s,\varepsilon)J_0^\prime(x),
\vspace{3mm} \\
\displaystyle
x=\frac{\tau}{\varepsilon},
\qquad
\tau=\int\limits_0^sP(t)dt.
\ea
\end{equation}
Here $u_i(s,\varepsilon)$, $w_i(s,\varepsilon)$ $(i=1,2)$ and
$P(t)$ are functions to be determined.

The advisability of representing solutions in the form
(\ref{trotsenko:2.10}) will be warranted if functions
$u_i(s,\varepsilon)$ and $w_i(s,\varepsilon)$ that we obtain are
regular and can be presented as a direct expansion in the
parameter $\varepsilon$.

Having substituted expressions (\ref{trotsenko:2.10}) into homogeneous
system of equations (\ref{trotsenko:2.8}) and set equal to zero
coef\/f\/icients of $J_0(x)$ and $J_0^\prime(x)$, we get four
equations for the functions
$u_i(s,\varepsilon)$, $w_i(s,\varepsilon)$ and $P(s)$:
\begin{equation} \label{trotsenko:2.11}
\ba{l}
\displaystyle
L_{11}(u_1)+L_{22}(w_1)+(\alpha_1Pu_2)'/\varepsilon+
\bigl(u_1P^2\alpha_1\bigr)/\varepsilon^2
\vspace{3mm} \\
\displaystyle
\qquad
+\alpha_1\left(u_2^\prime P-\frac{u_2P^2}{\tau}\right)\bigl/\varepsilon-
(\alpha_3w_2P)/\varepsilon=0,
\vspace{3mm} \\
\displaystyle
L_{11}(u_2)+L_{12}(w_2)+\left(\alpha_1Pu_2/\tau\right)'-
(\alpha_1Pu_1)'/\varepsilon+(\alpha_3Pw_1)/\varepsilon-
(\alpha_3Pw_2)/\tau
\vspace{3mm} \\
\displaystyle
\qquad
+\alpha_1u_2P^2\bigl(1/\varepsilon^2-1/\tau^2\bigr)+P\alpha_1
\bigl(u_1P/\tau-u_1^\prime\bigr)/\varepsilon+u_2^\prime
P\alpha_1/\tau=0, 
\vspace{3mm} \\
\displaystyle
\varepsilon^2\left[L_{21}(u_2)+L_{22}(w_2)+
\bigl(\beta_1w_2P/\tau\bigr)^\prime+
P\bigl(\beta_1w_2^\prime-\beta_3u_2\bigr)/\tau\right]+
rw_2r'+ru_2z'
\vspace{3mm} \\
\displaystyle
\qquad
+\varepsilon\left[\beta_3u_1P-(\beta_1w_1P)'+\beta_1P
\bigl(w_1P/\tau-w_1^\prime\bigr)\right]+\beta_1P^2w_2
\bigl(1-\varepsilon^2/\tau^2\bigr)=0,
\vspace{3mm} \\
\displaystyle
\varepsilon^2\bigl[L_{21}(u_1)+L_{22}(w_1)\bigr]+\varepsilon
\left[(\beta_1w_2P)'+\beta_1P\bigl(w_2^\prime-w_2P/\tau\bigr)-
\beta_3u_2P\right]
\vspace{3mm} \\
\displaystyle
\qquad
+w_1P^2\beta_1+rw_1r'+ru_1z'=0.
\ea
\end{equation}

We shall seek a solution to system of equations (\ref{trotsenko:2.11}) in
the form
\begin{equation} \label{trotsenko:2.12}
\ba{l}
\displaystyle
u_1=\sum\limits_{k=0}^\infty\varepsilon^{2k+2}u_1^{(k)}(s),
\qquad
u_2=\sum\limits_{k=0}^\infty\varepsilon^{2k+1}u_2^{(k)}(s),
\vspace{3mm} \\
\displaystyle
w_1=\sum\limits_{k=0}^\infty\varepsilon^{2k}w_1^{(k)}(s),
\qquad
w_2=\sum\limits_{k=0}^\infty\varepsilon^{2k+1}w_2^{(k)}(s).
\ea
\end{equation}

To determine functions $u_i^{(k)}$, $w_i^{(k)}$ and $P(s)$ we
substitute expansions (\ref{trotsenko:2.12}) into equations
(\ref{trotsenko:2.11})
and set  coef\/f\/icients of the same degrees of
$\varepsilon$ equal to zero.

In the zeroth-order approximation, i.e., for functions
$u_i^{(0)}$ and $w_i^{(0)}$, we  have
\begin{equation} \label{trotsenko:2.13}
\ba{l}
P^2\alpha_1u_1^{(0)}-\alpha_3 Pw_2^{(0)}=\psi,
\qquad \!
P\alpha_1u_2^{(0)}+\alpha_3w_1^{(0)}=0,
\qquad   \!
\bigl(P^2\beta_1+rr'\bigr)w_1^{(0)}=0,
\vspace{3mm} \\
\displaystyle
-\frac{d}{ds}\left(\beta_1Pw_1^{(0)}\right)+\beta_1P
\left(\frac{Pw_1^{(0)}}{\tau}-\frac{dw_1^{(0)}}{ds}\right)+
\bigl(P^2\beta_1+rr'\bigr)w_2^{(0)}+rz'u_2^{(0)}=0,
\vspace{3mm} \\
\displaystyle
\psi=-L_{12}\bigl(w_1^{(0)}\bigr)-
\frac{d}{ds}\left(\alpha_1Pu_2^{(0)}\right)-\alpha_1
\left(P\frac{du_2^{(0)}}{ds}-\frac1\tau u_2^{(0)}P^2\right).
\ea\!\!\!\!
\end{equation}
\begin{equation} \label{trotsenko:2.14}
u_2^{(0)}=-\frac{\alpha_3w_1^{(0)}}{\alpha_1P},
\qquad
P^2=-\frac{r}{\beta_1}\frac{dr}{ds}=-\frac{\lambda_1^2\cos\alpha}{T_1}.
\end{equation}
Here $\alpha$ is the angle formed by the outer normal to the
deformed shell and its symmetry axis.

It follows from (\ref{trotsenko:2.14}) that $P(s)$ is an even and pure
imaginary function for $\displaystyle
-\frac{\pi}{2}\leq\alpha\leq \frac{\pi}{2}$. The fourth equation
with regard for (\ref{trotsenko:2.13}) takes the form
\[ 
2\frac{dw_1^{(0)}}{ds}=\left(\frac P\tau-\frac{1}{\beta_1P}
\frac{d}{ds}(\beta_1P)+\varphi(s)\right)w_1^{(0)},
\qquad
\varphi(s)=\frac{\alpha_3}{\alpha_1}\frac{dz}{ds}\biggl/\frac{dr}{ds}.
\]

Hence up to an arbitrary constant, we obtain a regular
solution for function $w_1^{(0)}$:
\begin{equation} \label{trotsenko:2.15}
w_1^{(0)}=\sqrt{\frac{\tau}{\beta_1P}\exp\int \varphi(s)ds}.
\end{equation}

Therefore, after the zeroth-order approximation, the function
$w_2^{(0)}$ remains a free function that will be determined at the
second step of the asymptotic integration of equations~(\ref{trotsenko:2.11}).

Without to f\/inding explicit expressions for terms
in expansions (\ref{trotsenko:2.12}) of higher approximations,
one may conclude from the form of system (\ref{trotsenko:2.11})
with regard to expressions for its coef\/f\/icients,  that
$u_2^{(k)}(s)$ and $w_1^{(k)}(s)$ are even functions, and  
$u_1^{(k)}(s)$ and $w_2^{(k)}(s)$ are odd ones. Specif\/ic
expressions for these functions after the zeroth-order and f\/irst
approximations show that they are regular  for  $s=0$.

Putting together the results obtained, we conclude that a
solution to system (\ref{trotsenko:2.8}), bounded for $s=0$ and having the
large variability, can be presented in the form of the following formal
expansions: 
\begin{equation} \label{trotsenko:2.16}
\ba{l}
\displaystyle
u(s,\varepsilon)=I_0(y)\sum\limits_{k=0}^\infty
\varepsilon^{2k+2}u_1^{(k)}(s)+
I_0^\prime(y)\sum\limits_{k=0}^\infty
\varepsilon^{2k+1}u_2^{(k)}(s),
\vspace{3mm} \\
\displaystyle
w(s,\varepsilon)=I_0(y)\sum\limits_{k=0}^\infty
\varepsilon^{2k}w_1^{(k)}(s)+
I_0^\prime(y)\sum\limits_{k=0}^\infty
\varepsilon^{2k+1}w_2^{(k)}(s),
\ea
\end{equation}
where $I_0$ is a Bessel function of the order zero of pure imaginary
argument; $y$ is a function regular and odd in argument $s$.

The constructed asymptotic representations of solutions
(\ref{trotsenko:2.6}) and (\ref{trotsenko:2.16}) quantitatively present the
structure of a solution, that we use in a sequel to
construct systems of coordinate functions when solving the initial
problem by the method of variations.

Taking into account the asymptotics of expansions (\ref{trotsenko:2.16})
for $s\to 0$, and also the asymptotics of fucntions $I_n(x)$
for large values of the argument, after passing from normal and
tangent movements of the shell to its movements in axial and
radial directions, the general form of solutions to the initial
nonlinear problem can be presented as follows:
\begin{equation} \label{trotsenko:2.17}
\ba{l}
\displaystyle
z(s)=a_0+\varphi_p(s)\bigl(a_1+a_2s^2+a_3s^4+\cdots\bigr),
\vspace{3mm} \\
\displaystyle
r(s)=b_0s+\varphi_p(s)\bigl(b_1s+b_2s^3+b_3s^5+\cdots\bigr),
\vspace{3mm} \\
\displaystyle
\varphi_p(s)=I_0\left(\sum\limits_{k=0}^np_ks^{2k-1}\right),
\ea
\end{equation}
where $a_i$, $b_i$ and $p_i$ are arbitrary constants.

Subjecting expressions (\ref{trotsenko:2.17}) to boundary conditions
(\ref{trotsenko:1.5}), coordinate systems $\{u_k(s)\}$ and $\{v_k(s)\}$
take the form
\begin{equation} \label{trotsenko:2.18}
\ba{l}
\displaystyle
u_1(s,p_i)=\left(1-\frac{\varphi_p(s)}{\varphi_p(1)}\right),
\qquad
u_2(s,p_i)=\frac{\varphi_p(s)}{\varphi_p(1)}(s^2-1),
\vspace{3mm} \\
\displaystyle
u_k(s,p_i)=s^2u_{k-1}(s,p_i)
\qquad
(k=3,4,\ldots),
\qquad
v_k(s,p_i)=su_k(s,p_i).
\ea
\end{equation}

If constants $p_i$ are known, then solving nonlinear algebraic
system (\ref{trotsenko:1.9}) by means of the iteration Newton method, we
can construct an approximate solution for the problem. However,
the question of specif\/ic values of parameters $p_i$ for the
nonlinear boundary problem under consideration remains open.
Hence, from the conditions $\partial I/\partial p_i=0$,
$i=1,2,\ldots,n$ we can get  $n$ additional equations with respect to
parameters $p_i$. These equations are of the following form:
\begin{equation} \label{trotsenko:2.19}
\ba{l}
\displaystyle
\int\limits_0^1\left[ U(\lambda_1,\lambda_2)\frac{dz}{ds}
\frac{d}{ds}\left(\frac{dz}{dp_i}\right)-Q\lambda_2\frac{dr}{ds}
\frac{dz}{dp_i}+U(\lambda_1,\lambda_2)\frac{dr}{ds}\frac{d}{ds}
\left(\frac{dr}{dp_i}\right)\right.
\vspace{3mm} \\
\displaystyle
\qquad
\left. +\left(U(\lambda_2,\lambda_1)\frac{\lambda_2}{s}+Q\lambda_2
\frac{dz}{ds}\right)\frac{dr}{dp_i}\right] s\hspace{1.5pt}ds.
\ea
\end{equation}

It is convenient to look for a solution of system (\ref{trotsenko:2.19}) by
the generalized secant method~\cite{trotsenko:11:}. The whole algorithm is
reduced to that at each step of changing parameters $p_i$ (outer
iteration process) one needs to solve nonlinear algebraic system
(\ref{trotsenko:1.9}).

Therefore, unlike the traditional Rietz method of solving
nonlinear boundary problems, a peculiarity of the given approach
is that at the expense of choosing parameters $p_i$
that characterize the variability of required
solutions,  by using a
computer,  we construct a system of coordinate functions, that
approximates solutions to the initial problem in a way optimal in
some sense.

In conclusion, let us present results of calculations, giving an
idea of possibilities and ef\/fectiveness of the algorithm
proposed. To illustrate the convergence of successive
approximations, in Table 2 we write out  values of functions 
$z(s)$ and $r(s)$ and their f\/irst two derivatives depending on
the number $m$ of terms in expansions (\ref{trotsenko:1.7}) at the point
$s=0.9$. In the last column we give the relative error of an
approximate solution for problem (\ref{trotsenko:1.14}). Parameters of the
hydrostatic load in the coordinate system associated with the
center of the nondeformed membrane are $C=0.5$, $D=10$. Ratios
of constants in elastic potential~(\ref{trotsenko:1.2}) are chosen as
follows: 
\[
\Gamma_1=0.1,
\qquad
\Gamma_2=\Gamma_3=0.
\]
All calculations are performed with regard for only one
parameter $P_1$ in coordinate systems (\ref{trotsenko:2.18}).

\begin{center}
{\bf Table 2}\\[3mm]
\begin{tabular}{cccccccc}
\hline
\multicolumn{1}{c|}{$m$\rule[-7pt]{0pt}{20pt}} & 
\multicolumn{1}{c|}{$z\cdot 10$}
& \multicolumn{1}{c|}{$-z'$} & \multicolumn{1}{c|}{$-z''$} &
\multicolumn{1}{c|}{$r$} & \multicolumn{1}{c|}{$r'$} &
\multicolumn{1}{c|}{$-r''$} & \multicolumn{1}{c}{$\Delta$} \\
\hline
\rule{0pt}{13pt}%
1 & 0.32896 & 0.18816 & 1.9342 & 0.90532 & 0.97550 & 0.38008 &
$4\cdot 10^{-1}$ \\
2 & 0.36718 & 0.18415 & 2.4584 & 0.90694 & 0.99159 & 0.46372 &
$1\cdot 10^{-2}$ \\
3 & 0.36533 & 0.18145 & 2.3754 & 0.90692 & 0.99250 & 0.42182 &
$1\cdot 10^{-2}$ \\
4 & 0.36451 & 0.17825 & 2.3497 & 0.90693 & 0.99278 & 0.41433 &
$3\cdot 10^{-4}$ \\
5 & 0.36450 & 0.17830 & 2.3498 & 0.90693 & 0.99276 & 0.41456 &
$2\cdot 10^{-4}$ \\
6 & 0.36448 & 0.17841 & 2.3461 & 0.90693 & 0.99275 & 0.41404 &
$3\cdot 10^{-5}$ \\
\hline
\end{tabular}
\end{center}

Solutions possess also the analogous convergence at other points
from the interval of integration of the initial equations. Here
the high precision of an approximate solution is achieved at the expense
of both increasing the dimension of algebraic system
(\ref{trotsenko:1.9}) and choosing parameters of the optimal coordinate
system. It turns out that it is suf\/f\/icient to restrict oneself to
one parameter~$p_1$ only. Calculations show that requirements to the
precision of determining parameters $p_i$ concern with input
data of the hydrostatic load with high singularity in solutions
of the problem only. As the width of the boundary layer increases,
the rate of convergence of the algorithm is preserved, but inf\/luence
of the precision of determining parameters for the optimal
coordinate system on the f\/inal result decreases. There comes a
moment when for some value of the hydrostatic load we do not
need any more to organize the outer iteration process for
seeking for values of parameters $p_i$. In this case, the rate of
convergence of the given algorithm coincides with the rate of
convergence of the algorithm based on the use of an exponent
basis. It is worth noting that a polynomial basis for the
example considered above allows one to obtain the relative error
$\Delta$ for the f\/irst eight terms in expansions (\ref{trotsenko:1.7}), that
is equal to $\approx 10^{-2}$. The further increase of the
dimension of algebraic system (\ref{trotsenko:1.9}) in this case leads to
the loss of the stability of calculations.

While performing the calculating experiment, it came out that
functional (\ref{trotsenko:1.4}) in the space of parameters $p_i$ had only
one minimum. This means that for the approximation under
consideration there exists  only one optimal coordinate
system, which considerably simplify its search by means of a
computer. 

Therefore, the proposed approach allows one to obtain
approximate solutions with a suf\/fuciently high degree of
precision for the small dimension of an algebraic system both in
the presence of the singular perturbation of a problem and
without it.

\begin{figure}[t]

\centerline{\scalebox{1}{\includegraphics{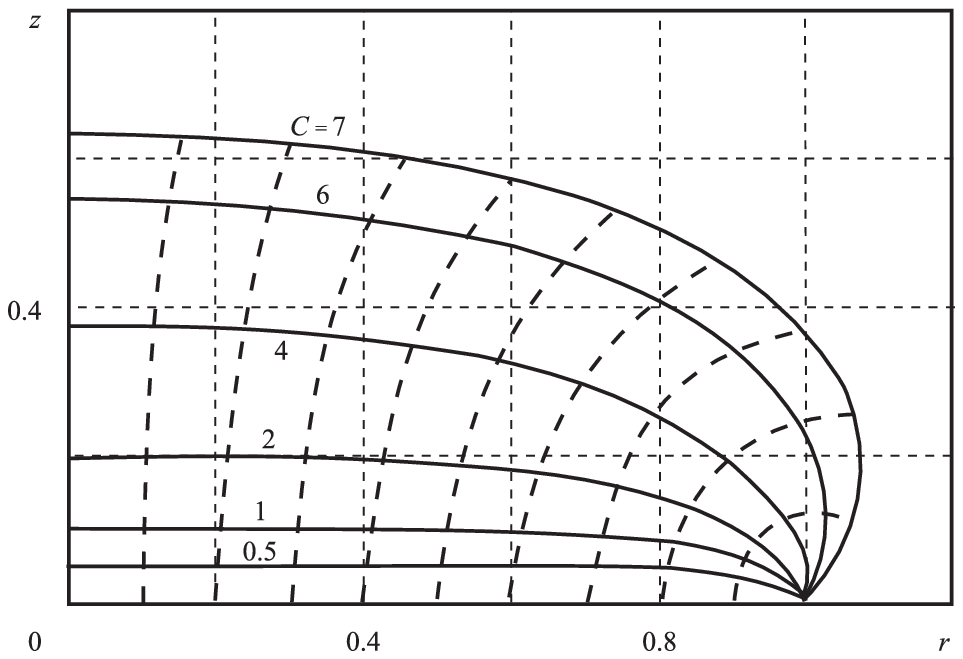}}}
\begin{center}
Figure 2
\end{center}

\vspace{-7mm}

\end{figure}

Meridional cuts of the deformed shell for dif\/ferent values of
the hydrostatic load $C$ are presented in Figure\,2. The other
input data coincide with those used for calculations presented
in Table~2. Dashed lines correspond to paths of passage of
points of the nondeformed membrane to points of its deformed
surface. Curves of dependence of principal degrees of
lengthenings $\lambda_1$ and stresses $T_1$ (continuous curves)
and also $\lambda_2$ and stresses~$T_2$ (dashed curves) upon the
parameter $s$ for dif\/ferent values of the parameter $C$ are
presented in Figures\,3 and 4. From f\/igures one can see that at the
starting stage of deforming the membrane, the main part of its
surface is nearly plane that is in the homogeneous biaxial
strained state. Only in a neighbourhood of the support contour of
the shell its surface and intererior stresses undergo
considerable changes. As the constant component of the
hydrostatic load $C$ increases, the plate part of the surface
gradually decreases.

\begin{figure}[t]
\centerline{\scalebox{1}{\includegraphics{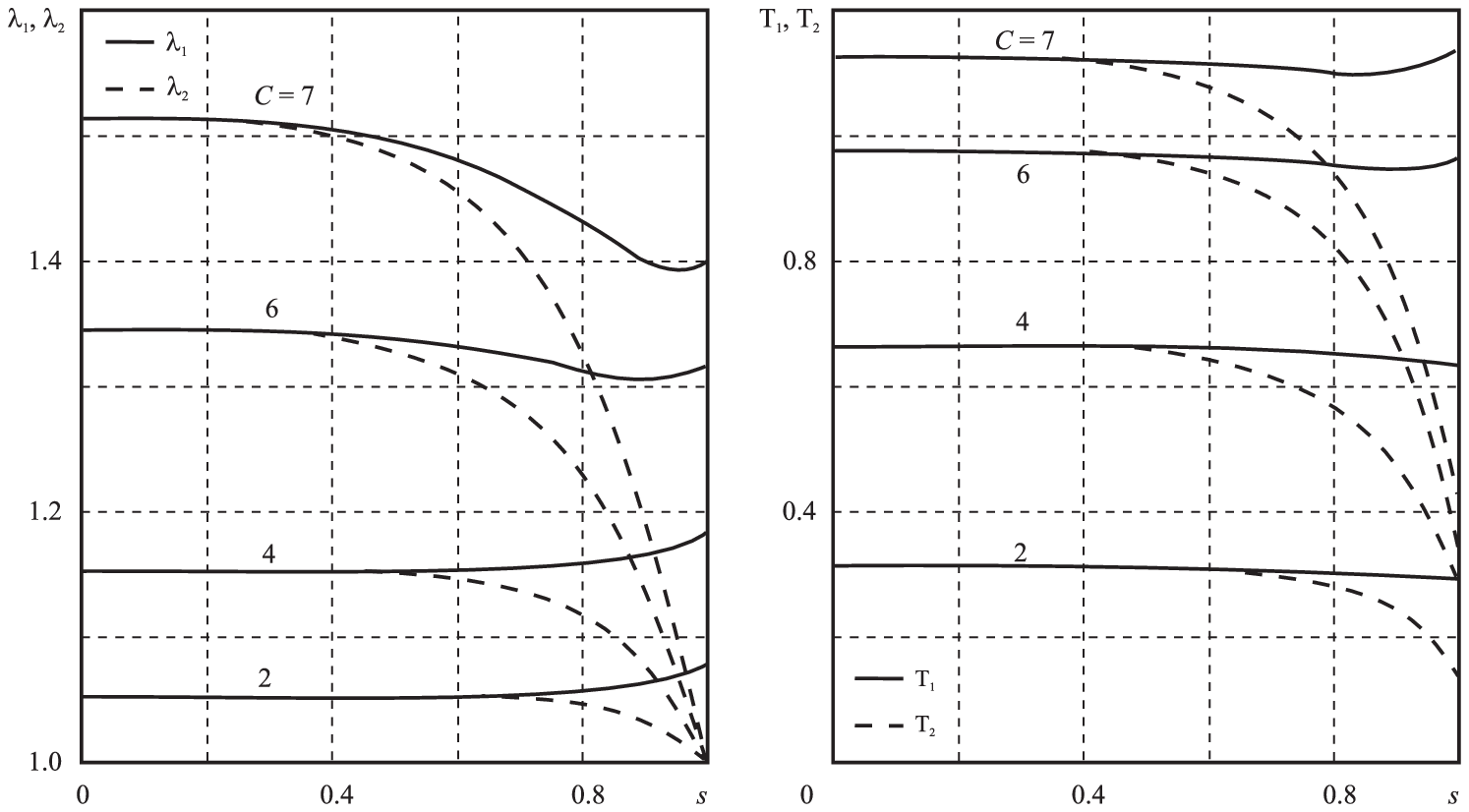}}}
\begin{center}
Figure 3     \hspace{65mm} Figure 4
\end{center}

\vspace{-7mm}

\end{figure}

Since it is convenient to investigate behaviour of material of
elastomers for large deformations by means of experiments on
homogeneous biaxial extension of a thin sheet~\cite{trotsenko:1:} (that
often causes some dif\/f\/iculties when realizing such a type of
deformations for experimental tests), the presence of the homogeneous
strained state in the central part of the deformed membrane
under the load of heavy liquid may be used in establishing a
functional dependence of an elastic potential upon deformation invariants.

\subsection*{Acknowledgements}

The author  is grateful to the  DFG for partial f\/inancial support.

\label{trotsenko-lp}

\end{document}